\documentclass[a4paper,12pt]{article}
\usepackage[T2A]{fontenc}
\usepackage[english]{babel}
\usepackage{graphicx}
\usepackage{amsmath}
\usepackage{mathrsfs}
\usepackage{amssymb}

\begin{document}

\begin{center}

{\bf Stability and convergence of difference schemes approximating a
two-parameter nonlocal boundary value problem for time-fractional
diffusion equation}

\vspace{2mm}

{ Anatoly A. Alikhanov}

 \vspace{2mm}

 {Kabardino-Balkarian State
University, ul. Chernyshevskogo 173,
 Nalchik,  360004,   Russia, e-mail: aaalikhanov@gmail.com}

\end{center}

\begin{abstract}
Difference schemes for the time-fractional diffusion equation with
variable coefficients and nonlocal boundary conditions containing
real parameters $\alpha$ and $\beta$ are considered. By the method
of energy inequalities, for the solution of the difference problem,
we obtain a priori estimates, which imply the stability and
convergence of these difference schemes.
\end{abstract}


\section{Introduction}
Consider the nonlocal boundary value problem

\begin{equation}\label{ur1}
\partial_{0t}^{\gamma}u=\frac{\partial }{\partial
x}\left(k(x)\frac{\partial u}{\partial x}\right)+f(x,t) ,\quad
0<x<1,\quad 0<t\leq T,
\end{equation}

\begin{equation}
u(0,t)=\alpha u(1,t),\quad
 k(1)u_x(1,t)=\beta
k(0)u_x(0,t)+\mu(t),\quad 0\leq t\leq T, \label{ur2}
\end{equation}

\begin{equation}
u(x,0)=u_0(x),\quad 0\leq x\leq 1, \label{ur3}
\end{equation}
where $0<c_1\leq k(x)\leq c_2$ for all $x\in [0,1]$, $\alpha$ and
$\beta$ are real numbers such that $\alpha\beta>0$, and $\mu(t)$
$\in C[0,T]$, $\partial_{0t}^{\gamma}
u(x,t)=\int_0^tu_{\tau}(x,\tau)(t-\tau)^{-\gamma}d\tau/\Gamma(1-\gamma)$
 is a  Caputo fractional derivative of order $\gamma$, $0<\gamma<1$ \cite{Nakh:03, Cap:69}.

The first nonlocal condition in (\ref{ur2}) can be replaced by the
inhomogeneous condition $u(0,t)=\alpha u(1,t)+\mu_1(t)$, however, if
$\mu_1(t) \in C^1[0,T]$, then the simple change of variables
$u(x,t)=v(x,t)+(1-x)\mu_1(t)$ reduces this problem to the considered
one.

The existence of the solution for the initial boundary value problem
of a number of fractional order differential equations has been
proved in \cite{Luchko_ibv, Luchko_mt}.

We introduce the space grid $\omega_h=\{x_i=ih\}_{i=0}^{N}$, and the
time grid $\omega_\tau=\{t_n=n\tau\}_{n=0}^{Nt}$ with increments
$h=1/N$ and $\tau=T/Nt$. Set
$$a_i=k(x_i-0.5h),\quad \varphi_i^n=f(x_i,t_n+\sigma\tau), \quad
y_i^n=y(x_i,t_n), \quad y_{\bar x,i}^n=(y_i^n-y_{i-1}^n)/h,
$$
$$y_{ x,i}^n=(y_{i+1}^n-y_{i}^n)/h, \quad (ay_{\bar
x})_{x,i}^n=(a_{i+1}y_{i+1}^n-(a_{i+1}+a_i)y_i^n+a_iy_{i-1}^n)/h^2,
$$
$$y_{t,i}^n=(y_i^{n+1}-y_i^n)/\tau, \quad v_i^{(\sigma)}=\sigma
v_i^{n+1}+(1-\sigma)v_i^n, \quad 0\leq \sigma\leq 1,
$$
$\Delta_{0t_{n+1}}^{\gamma}
y_i=\sum_{s=0}^n(t_{n-s+1}^{1-{\gamma}}-t_{n-s}^{1-{\gamma}})y_{
t,i}^s/\Gamma(2-{\gamma})$ is a difference analogue of Caputo
fractional derivative of order $\gamma$, $0<\gamma<1$
\cite{ShkhTau:06}.

It has been shown that if the function $v(t)\in C^2[0,T]$ then
$\partial_{0t_{j+1}}^{\alpha}v=\Delta_{0t_{j+1}}^{\alpha}v+O(\tau)$
\cite{ShkhTau:06}. This result can be improved if $v(t)\in
C^3[0,T]$.

{\bf Lemma 1.} For any function $v(t)\in C^3[0,T]$ the following
equality takes place
\begin{equation}
\partial_{0t_{j+1}}^{\alpha}v=\Delta_{0t_{j+1}}^{\alpha}v+O(\tau^{2-\alpha}),
\quad 0<\alpha<1. \label{ur3.5}
\end{equation}

{\bf Proof.}  For $t=t_{j+1}$, $j=0,1,\ldots,j_0-1$ one has
$$
\partial_{0t_{j+1}}^{\alpha}v=\frac{1}{\Gamma(1-\alpha)}\int\limits_{0}^{t_{j+1}}\frac{v'(\eta)d\eta}{(t_{j+1}-\eta)^{\alpha}}=
\frac{1}{\Gamma(1-\alpha)}\sum\limits_{s=0}^{j}\int\limits_{t_s}^{t_{s+1}}\frac{v'(\eta)d\eta}{(t_{j+1}-\eta)^{\alpha}}=
$$
$$
=\frac{1}{\Gamma(1-\alpha)}\sum\limits_{s=0}^{j}\int\limits_{t_s}^{t_{s+1}}\frac{v'(t_{s+1/2})+
v''(t_{s+1/2})(\eta-t_{s+1/2})+O((\eta-t_{s+1/2})^2)}{(t_{j+1}-\eta)^{\alpha}}d\eta=
$$
$$
=\Delta_{0t_{j+1}}^{\alpha}v+\frac{1}{\Gamma(1-\alpha)}\sum\limits_{s=0}^{j}v''(t_{s+1/2})\int\limits_{t_s}^{t_{s+1}}\frac{
\eta-t_{s+1/2}}{(t_{j+1}-\eta)^{\alpha}}d\eta+O(\tau^2).
$$
Let us estimate the value
$$
\left|\frac{1}{\Gamma(1-\alpha)}\sum\limits_{s=0}^{j}v''(t_{s+1/2})\int\limits_{t_s}^{t_{s+1}}\frac{
\eta-t_{s+1/2}}{(t_{j+1}-\eta)^{\alpha}}d\eta\right|\leq
\frac{M}{\Gamma(1-\alpha)}\sum\limits_{s=0}^{j}\left|\int\limits_{t_s}^{t_{s+1}}\frac{
\eta-t_{s+1/2}}{(t_{j+1}-\eta)^{\alpha}}d\eta\right|=
$$
$$
=\frac{M}{\Gamma(1-\alpha)}\sum\limits_{s=0}^{j}\left|\int\limits_{t_{s+1/2}}^{t_{s+1}}\frac{
\eta-t_{s+1/2}}{(t_{j+1}-\eta)^{\alpha}}d\eta-\int\limits_{t_s}^{t_{s+1/2}}\frac{
t_{s+1/2}-\eta}{(t_{j+1}-\eta)^{\alpha}}d\eta\right|=
$$
$$
=\frac{2^{\alpha}M\tau^{2-\alpha}}{4\Gamma(1-\alpha)}\sum\limits_{s=0}^{j}\left(\int\limits_{0}^{1}\frac{
zdz}{(2(j-s)+1-z)^{\alpha}}-\int\limits_{0}^{1}\frac{
zdz}{(2(j-s)+1+z)^{\alpha}}\right)=
$$
$$
=\frac{2^{\alpha}M\tau^{2-\alpha}}{4\Gamma(1-\alpha)}\int\limits_{0}^{1}z\sum\limits_{s=0}^{j}\left(\frac{
1}{(2s+1-z)^{\alpha}}-\frac{1}{(2s+1+z)^{\alpha}}\right)dz=
$$
$$
=\frac{2^{\alpha}M\tau^{2-\alpha}}{4\Gamma(1-\alpha)}\int\limits_{0}^{1}\left(\frac{
z}{(1-z)^{\alpha}}-\frac{z}{(2j+1+z)^{\alpha}}\right)dz-
$$
$$
-\frac{2^{\alpha}M\tau^{2-\alpha}}{4\Gamma(1-\alpha)}\int\limits_{0}^{1}z\sum\limits_{s=1}^{j}\left(\frac{
1}{(2s-1+z)^{\alpha}}-\frac{1}{(2s+1-z)^{\alpha}}\right)dz\leq
$$
$$
\leq\frac{2^{\alpha}M\tau^{2-\alpha}}{4\Gamma(1-\alpha)}\int\limits_{0}^{1}\frac{
zdz}{(1-z)^{\alpha}}=\frac{2^{\alpha}M}{4\Gamma(3-\alpha)}\tau^{2-\alpha},
$$
where $M=\max\limits_{0\leq t\leq T}|v''(t)|$. The proof of Lemma 1
is complete.

Consider the weighted scheme

\begin{equation}
\Delta_{0t_{n+1}}^{\gamma} y_i-(ay_{\bar
x}^{(\sigma)})_{x,i}=\varphi_i^n,\quad i=1,2,...,N-1, \label{ur4}
\end{equation}

\begin{equation}
\begin{cases}
y_0^{n+1}-\alpha y_N^{n+1}=0, \\
\beta \Delta_{0t_{n+1}}^{\gamma} y_0+\Delta_{0t_{n+1}}^{\gamma}
y_N+\dfrac{2}{h}\left(a_Ny_{\bar x,
N}^{(\sigma)}-\beta a_1 y_{ x, 0}^{(\sigma)}\right)=\dfrac{2}{h}\mu(t_{n+\sigma})+\varphi_N+\beta\varphi_0,\\
\end{cases}
\label{ur5}
\end{equation}

\begin{equation}
y_i^0=u_0(x_i).
 \label{ur6}
\end{equation}
The difference scheme (\ref{ur4})--(\ref{ur6}) has approximation
order $O(\tau^{m_\sigma}+h^2)$ where $m_\sigma=1$ if $0\leq\sigma<1$
and $m_\sigma=2-\alpha$ if $\sigma=1$ \cite{Samar:77}.

The nonlocal boundary value problem with the boundary conditions
$u(b,t)=\rho u(a,t)$, $u_x(b,t)=\sigma u_x(a,t)+\tau u(a,t)$ for the
simplest equations of mathematical physics, referred to as
conditions of the second class, was studied in the monograph
\cite{Steclov}. Results in the case in which $ \rho\sigma-1=0$ and
$\rho\tau\leq 0$ were obtained there. Difference schemes for problem
(\ref{ur1})--(\ref{ur3}) with $\alpha=\beta$ and $\gamma=1$ (the
classical diffusion equation) were studied in \cite{Gul1}. In this
case, the operator occurring in the elliptic part is self-adjoint.
Self-adjointness permits one to use general theorems on the
stability of two-layer difference schemes in energy spaces and
consider difference schemes for equations with variable
coefficients. Stability criteria for difference schemes for the heat
equation with nonlocal boundary conditions were studied in
\cite{Gul1.5, Gul1.6, Gul1.7, Gul1.8, Gul1.9}. The difference
schemes considered in these papers have the specific feature that
the corresponding difference operators are not self-adjoint. The
method of energy inequalities was developed in \cite{Alikh:08,
Alikh:10_2, Alikh:13} for the derivation of a priory estimates for
solutions of difference schemes for the classical diffusion equation
with variable coefficients in the case of nonlocal boundary
conditions. Using the energy inequality method, a priory estimates
for the solution of the Dirichlet and Robin boundary value problems
for the fractional, variable and distributed order diffusion
equation with Caputo fractional derivative have been obtained
\cite{Alikh:10, Alikh:12, Alikh:13_3}. A priori estimates for the
difference problems analyzed in \cite{ShkhTau:06, ShkhLaf:09,
ShkhLaf:10} by using the maximum principle imply the stability and
convergence of these difference schemes.

The existence and uniqueness of solutions for fractional ordinary
differential equations with various kinds of the fractional
derivative and nonlocal boundary conditions have been proven
\cite{odu1, odu2, odu3}.

Numerical methods for solving fractional diffusion equations with
classical boundary value problems and various kinds of the
fractional order derivative have been proposed \cite{Liu:09,
Liu:09_2, Liu:10, Liu:12_1, Liu:12_2}.

\section{A priori estimate for the differential problem.}

{\bf Lemma 1.} \cite{Alikh:12} For any function $v(t)$ absolutely
continuous on $[0,T]$ the following equality takes place:

\begin{equation}
v(t)\partial_{0t}^{\nu} v(t)=\frac{1}{2}\partial_{0t}^{\nu}v^2(t)+
\frac{\nu}{2\Gamma(1-\nu)}\int\limits_{0}^{t}\frac{d\xi}{(t-\xi)^{1-\nu}}
\left(\int\limits_{0}^{\xi}\frac{v'(\eta)d\eta}{(t-\eta)^\nu}\right)^2,
 \label{ur5555}
\end{equation}
where $0<\nu<1$.

{\bf Theorem 1.} The solution of the nonlocal boundary value problem
(\ref{ur1})--(\ref{ur3}) satisfies the identity

$$
\frac{1}{2}\partial_{0t}^{\gamma}\int\limits_{0}^{1}(1+\delta
p(x))u^2(x,t)dx+\frac{\gamma
}{2\Gamma(1-\gamma)}\int\limits_{0}^{1}\left(1+\delta
p(x)\right)dx\int\limits_{0}^{t}
\frac{\left(\int_{0}^{\xi}\frac{\frac{\partial u}{\partial
\eta}(x,\eta)d\eta}{(t-\eta)^{\gamma}}\right)^2d\xi}{(t-\xi)^{1-\gamma}}+
$$
$$
+\int\limits_{0}^{1}(1+\delta
p(x))k(x)u_x^2(x,t)dx+\frac{\delta}{2}(\alpha^2-1) u^2(1,t)=
$$
\begin{equation}
=\int\limits_{0}^{1}(1+\delta p(x))u(x,t)f(x,t)dx+
u(1,t)\mu(t),\label{ur7}
\end{equation}
where $p(x)=\int\limits_{x}^{1}k^{-1}(s)ds$,
$\delta=(\beta\alpha^{-1}-1)/p(0)$.

{\bf Proof.} Let us multiply (\ref{ur1}) by $u(x,t)$ and integrate
the resulting relation over $x$ from $0$ to $1$:

\begin{equation}
\int\limits_{0}^{1}u(x,t)\partial_{0t}^{\gamma}u(x,t)dx-\int\limits_{0}^{1}(k(x,t)u_x(x,t))_xu(x,t)dx=
\int\limits_{0}^{1}u(x,t)f(x,t)dx. \label{ur8}
\end{equation}

This, together with the nonlocal boundary conditions (\ref{ur2}) and
equality (\ref{ur5555}), implies the relation

$$
\frac{1}{2}\partial_{0t}^{\gamma}\int\limits_{0}^{1}u^2(x,t)dx+\frac{\gamma}{2\Gamma(1-\gamma)}\int\limits_{0}^{1}dx\int\limits_{0}^{t}
\frac{d\xi}{(t-\xi)^{1-\gamma}}\left(\int\limits_{0}^{\xi}\frac{\frac{\partial
u}{\partial \eta}(x,\eta)d\eta}{(t-\eta)^{\gamma}}\right)^2+
$$
\begin{equation}
+\int\limits_{0}^{1}k(x)u_x^2(x,t)dx=\int\limits_{0}^{1}u(x,t)f(x,t)dx+\left(\frac{\beta}{\alpha}-1\right)k(0)u_x(0,t)u(0,t)
+\frac{1}{\alpha}u(0,t)\mu(t). \label{ur9}
\end{equation}

We multiply (\ref{ur1}) by $u(x,t)$ and integrate with respect to
$s$ from $0$ to $x$,

\begin{equation}
\int\limits_{0}^{x}u(s,t)\partial_{0t}^{\gamma}u(s,t)ds-\int\limits_{0}^{x}(k(s)u_s(s,t))_su(s,t)ds=
\int\limits_{0}^{x}u(s,t)f(s,t)ds. \label{ur10}
\end{equation}

Hence we obtain the identity

$$
\frac{1}{2}\partial_{0t}^{\gamma}\int\limits_{0}^{x}u^2(s,t)ds+
\frac{\gamma}{2\Gamma(1-\gamma)}\int\limits_{0}^{x}ds\int\limits_{0}^{t}
\frac{d\xi}{(t-\xi)^{1-\gamma}}\left(\int\limits_{0}^{\xi}\frac{\frac{\partial
u}{\partial \eta}(s,\eta)d\eta}{(t-\eta)^{\gamma}}\right)^2+
$$
\begin{equation}
+\int\limits_{0}^{x}k(s)u_s^2(s,t)ds=\int\limits_{0}^{x}u(s,t)f(s,t)ds+k(x)u_x(x,t)u(x,t)-k(0)u_x(0,t)u(0,t).
\label{ur11}
\end{equation}

We divide (\ref{ur11}) by $k(x)$ and integrate with respect to $x$
from $0$ to $1$,

$$
\frac{1}{2}\partial_{0t}^{\gamma}\int\limits_{0}^{1}p(x)u^2(x,t)dx+
\frac{\gamma}{2\Gamma(1-\gamma)}\int\limits_{0}^{1}p(x)dx\int\limits_{0}^{t}
\frac{d\xi}{(t-\xi)^{1-\gamma}}\left(\int\limits_{0}^{\xi}\frac{\frac{\partial
u}{\partial \eta}(x,\eta)d\eta}{(t-\eta)^{\gamma}}\right)^2+
$$
$$
+\int\limits_{0}^{1}p(x)k(x)u_x^2(x,t)dx=\int\limits_{0}^{1}p(x)u(x,t)f(x,t)dx+
$$
\begin{equation}
+\frac{1}{2}(1-\alpha^2)u^2(1,t)-p(0)k(0)u_x(0,t)u(0,t).
\label{ur12}
\end{equation}

By multiplying identity (\ref{ur12}) by
$\delta=(\beta\alpha^{-1}-1)/p(0)$ and by adding relation
(\ref{ur9}), we obtain identity (\ref{ur7}). The proof of Theorem 1
is complete.

{\bf Theorem 2.} If the condition
$(\beta\alpha^{-1}-1)(\alpha^2-1)\geq0$ is satisfied, then the
solution of problem(\ref{ur1})--(\ref{ur3}) with $f(x,t)\equiv 0$
and $\mu(t)\equiv 0$ satisfies the estimate
\begin{equation}
\|u^2(x,t)\|_0^2+2c_1D_{0t}^{-\gamma}\|u_x(x,s)\|_0^2
\leq\max\left\{\frac{\alpha}{\beta},\frac{\beta}{\alpha}\right\}\|u_0(x)\|_0^2,
\label{ur13}
\end{equation}
where $\|u(x,t)\|_0^2=\int\limits_{0}^{1}u^2(x,t)dx$,
$D_{0t}^{-\nu}u(x,t)=\int\limits_{0}^{t}(t-s)^{\nu-1}u(x,s)ds/\Gamma(\nu)$
is the fractional Riemann--Liouville integral of order  $\nu>0$.

{\bf Proof.} By applying the fractional integration operator
$D_{0t}^{-\gamma}$ to both sides of equality (\ref{ur7})  and by
using the inequalities $\min\{1,\beta\alpha^{-1}\}\leq 1+\delta
p(x)\leq \max\{1,\beta\alpha^{-1}\}$,  we obtain the a priori
estimate (\ref{ur13}). The proof of Theorem 2 is complete.

The inequality $(\beta\alpha^{-1}-1)(\alpha^2-1)\geq0$ is equivalent
to the system of inequalities $|\beta|\geq |\alpha|\geq1$ or
$|\beta|\leq |\alpha|\leq1$. Note that one can obtain a priori
estimates for $|\alpha|\geq |\beta|\geq1$ or $|\alpha|\leq
|\beta|\leq1$ by multiplying equation (\ref{ur1}) by $(k(x)u_x)_x$.
To avoid similar calculations, we assume sufficient smoothness of
the solution and the input data and show that this holds for the
differential problem (\ref{ur1})--(\ref{ur3}). Set
$w(x,t)=k(x)u_x(x,t)$. Then the function $w(x,t)$ satisfies the
problem

\begin{equation}
\partial_{0t}^{\gamma}w=k(x)w_{xx}+k(x)f_x(x,t),
 \label{ur16.1}
\end{equation}
\begin{equation}
w(0,t)=\frac{1}{\beta}w(1,t)-\frac{1}{\beta}\mu(t),\quad
w_x(1,t)=\frac{1}{\alpha}w_x(0,t)+\frac{1}{\alpha}f(0,t)-f(1,t),
\label{ur16.2}
\end{equation}
\begin{equation}
w(x,0)=k(0)u_0'(x). \label{ur16.3}
\end{equation}

Obviously, for problem (\ref{ur16.1})--(\ref{ur16.3}) with
$f(x,t)\equiv 0$ and $\mu(t)\equiv 0$, one can readily obtain an a
priori estimate of the form (\ref{ur13}) for $|\alpha|\geq
|\beta|\geq1$ or $|\alpha|\leq |\beta|\leq1$.

We have the inequalities
\begin{equation}
u(1,t)\mu(t)\leq
\frac{\varepsilon}{2}u^2(1,t)+\frac{1}{2\varepsilon}\mu^2(t),
\label{ur13.5}
\end{equation}

$$
\int\limits_{0}^{1}(1+\delta
p(x))u(x,t)f(x,t)dx\leq\frac{\varepsilon}{2}u^2(1,t)+
$$
\begin{equation}
+\frac{\varepsilon_1}{2}\int\limits_{0}^{1}(1+\delta
p(x))k(x)u_x^2(x,t)dx+\left(\frac{\gamma_1}{2\varepsilon}+\frac{\gamma_1\gamma_2}{2\varepsilon_1}\right)\int\limits_{0}^{1}(1+\delta
p(x))f^2(x,t)dx, \label{ur14}
\end{equation}
where $\varepsilon$, $\varepsilon_1>0$,
$\gamma_1=\int\limits_{0}^{1}(1+\delta p(x))dx$,
$\gamma_2=\int\limits_{0}^{1}(1+\delta p(x))^{-1}k^{-1}(x)dx$.

Inequality (\ref{ur13.5}) is obvious. Let us prove inequality
(\ref{ur14}). By virtue of the relation $
u(x,t)=u(1,t)-\int\limits_{x}^{1}u_s(s,t)ds$, we have
$$
\int\limits_{0}^{1}(1+\delta
p(x))u(x,t)f(x,t)dx=\int\limits_{0}^{1}(1+\delta
p(x))f(x,t)\left(u(1,t)-\int\limits_{x}^{1}u_s(s,t)ds\right)dx=
$$
$$
=u(1,t)\int\limits_{0}^{1}(1+\delta
p(x))f(x,t)dx-\int\limits_{0}^{1}u_x(x,t)dx\int\limits_{0}^{x}(1+\delta
p(s))f(s,t)ds\leq
$$
$$
\leq\frac{\varepsilon}{2}u^2(1,t)+\frac{\gamma_1}{2\varepsilon}\int\limits_{0}^{1}(1+\delta
p(x))f^2(x,t)dx+
$$
$$
+\int\limits_{0}^{1}\sqrt{(1+\delta
p(x))k(x)}|u_x(x,t)|dx\int\limits_{0}^{x}\frac{(1+\delta
p(s))}{\sqrt{(1+\delta p(x))k(x)}}|f(s,t)|ds\leq
$$
$$
\leq\frac{\varepsilon}{2}u^2(1,t)+\frac{\gamma_1}{2\varepsilon}\int\limits_{0}^{1}(1+\delta
p(x))f^2(x,t)dx+
$$
$$
+\frac{\varepsilon_1}{2}\int\limits_{0}^{1}(1+\delta
p(x))k(x)u_x^2(x,t)dx+\frac{\gamma_1\gamma_2}{2\varepsilon_1}\int\limits_{0}^{1}(1+\delta
p(x))f^2(x,t)dx.
$$
The proof of inequality (\ref{ur14}) is complete.

{\bf Theorem 3.} If the condition
$(\beta\alpha^{-1}-1)(\alpha^2-1)>0$ is satisfied, then the solution
of problem (\ref{ur1})--(\ref{ur3}) satisfies the a priori estimates
\begin{equation}
\|u^2(x,t)\|_0^2+D_{0t}^{-\gamma}\|u_x(x,t)\|_0^2 \leq
M\left(D_{0t}^{-\gamma}\|f(x,t)\|_0^2+D_{0t}^{-\gamma}\mu^2(t)+\|u_0(x)\|_0^2\right),
\label{ur15}
\end{equation}
where $M=M(\alpha,\beta,c_1)>0$ is a known constant independent of
$t$.

{\bf Proof.} Identity (\ref{ur7}), together with inequalities
(\ref{ur13.5}) end (\ref{ur14}) for
$\varepsilon=(\beta\alpha^{-1}-1)(\alpha^2-1)$ and
$\varepsilon_1=1$, implies the inequality

$$
\frac{1}{2}\partial_{0t}^{\gamma}\int\limits_{0}^{1}(1+\delta
p(x))u^2(x,t)dx+\frac{\gamma}{\Gamma(1-\gamma)}\int\limits_{0}^{1}(1+\delta
p(x))dx\int\limits_{0}^{t}
\frac{\left(\int_{0}^{\xi}\frac{\frac{\partial u}{\partial
\eta}(x,\eta)d\eta}{(t-\eta)^{\gamma}}\right)^2d\xi}{(t-\xi)^{1-\gamma}}+
$$
\begin{equation}
+\frac{c_1}{2}\int\limits_{0}^{1}(1+\delta p(x))u_x^2(x,t)dx\leq
M_1\left(\int\limits_{0}^{1}(1+\delta
p(x))f^2(x,t)dx+\mu^2(t)\right),\label{ur16}
\end{equation}
where $M_1=M_1(\alpha,\beta,c_1)>0$ is a known constant independent
of $t$.

By applying the fractional integration operator $D_{0t}^{-\gamma}$
to both sides of inequality (\ref{ur16}) and by taking into account
the inequalities $\min\{1,\beta\alpha^{-1}\}\leq 1+\delta p(x)\leq
\max\{1,\beta\alpha^{-1}\}$,  we obtain the a priori estimate
(\ref{ur15}). The proof of Theorem 3 is complete.

\subsection{A priori estimates for the difference problem.}

{\bf Lemma 2.} \cite{Alikh:12} For any function $y(t)$ defined on
the grid $\bar\omega_{\tau}$ one has the equalities

\begin{equation}\label{ur32.5555}
 y^{n+1}\Delta_{0t}^\nu y = \frac{1}{2}\Delta_{0t}^\nu (y^2)
 +\frac{\tau^\nu \Gamma(2-\nu)}{2}(\Delta_{0t}^\nu y)^2+J_1(y),
\end{equation}

\begin{equation}\label{ur32.6666}
 y^{n}\Delta_{0t}^\nu y = \frac{1}{2}\Delta_{0t}^\nu (y^2)
 -\frac{\tau^\nu \Gamma(2-\nu)}{2(2-2^{1-\nu})}(\Delta_{0t}^\nu
 y)^2+J_2(y),
\end{equation}
where
$$
J_1(y)= \frac{1}{2\Gamma(2-\nu)}\sum\limits_{k=0}^{n-1}
\tau\left((t_{n-k+1}^{1-\nu}-t_{n-k}^{1-\nu})^{-1}
-(t_{n-k}^{1-\nu}-t_{n-k-1}^{1-\nu})^{-1}\right)\left(\zeta^{k+1}\right)^2,
$$
$$
J_2(y)=\frac{\tau^{\nu}(2^{1-\nu}-1)}{2\Gamma(2-\nu)(2-2^{1-\nu})}\left(\zeta^{n+1}
-\frac{2-2^{1-\nu}}{2^{1-\nu}-1}\zeta^{n}\right)^2+
$$
$$
+\frac{1}{2\Gamma(2-\nu)}\sum\limits_{k=0}^{n-2}
\tau\left((t_{n-k+1}^{1-\nu}-t_{n-k}^{1-\nu})^{-1}
-(t_{n-k}^{1-\nu}-t_{n-k-1}^{1-\nu})^{-1}\right)\left(\zeta^{k+1}\right)^2,
$$
$\zeta^{k+1}=\sum_{s=0}^k(t_{n-s+1}^{1-\nu}-t_{n-s}^{1-\nu})y_{
t}^s$, $ J_1(y)\geq0$, $J_2(y)\geq0 $, $0<\nu<1$. Here we consider
the sums to be equal to zero if the upper summation limit is less
than the lower one.

{\bf Lemma 3.} For any nonnegative function $v(x)$ defined on the
grid $\bar \omega_h$ and any solution $y(x,t)$ of equation
(\ref{ur4}) with zero right--hand side $\varphi\equiv 0$, one has
the inequality
$$
\|\sqrt{v}\sqrt{a}y_{\bar
x}^{(\sigma)}]|_0^2\geq\frac{h^2}{4c_2}\|\sqrt{v}\Delta_{0t}^{\gamma}
y\|_0^2+
$$
\begin{equation}
+\frac{h}{2c_2}\sum\limits_{i=1}^{N}v_{\bar x,i}(a_iy_{\bar
x,i}^{(\sigma)})^2h+\frac{h}{2c_2}\left(v_N(a_Ny_{\bar
x,N}^{(\sigma)})^2+v_0(a_1y_{x,0}^{(\sigma)})^2\right),
 \label{ur32}
\end{equation}
where $\|y\|_0^2=\sum_{i=1}^{N-1}y_i^2h$,
$\|y]|_0^2=\sum_{i=1}^{N}y_i^2h$.

{\bf Proof.} Since $y(x,t)$ is a solution of equation (\ref{ur4})
with zero right--hand side, it follows that, for each nonnegative
function $v(x)$ one has the relation
$\sqrt{v_i}\Delta_{0t}^{\gamma_i} y_i=\sqrt{v_i}(ay_{\bar
x}^{(\sigma)})_{x,i}$ for all $i=1,2,\ldots,N-1$; consequently,

$$
\|\sqrt{v}\Delta_{0t}^{\gamma} y\|_0^2=\|\sqrt{v}(ay_{\bar
x}^{(\sigma)})_x\|_0^2=
\frac{1}{h^2}\sum\limits_{i=1}^{N-1}v_i\left(a_{i+1}y_{\bar
x,i+1}^{(\sigma)}-a_{i}y_{\bar x,i}^{(\sigma)}\right)^2h\leq
$$
$$
\leq\frac{2}{h^2}\sum\limits_{i=1}^{N-1}v_i\left((a_{i+1}y_{\bar
x,i+1}^{(\sigma)})^2+(a_{i}y_{\bar x,i}^{(\sigma)})^2\right)h
=\frac{2}{h^2}\sum\limits_{i=1}^{N-1}\left(v_{i+1}(a_{i+1}y_{\bar
x,i+1}^{(\sigma)})^2+v_{i}(a_{i}y_{\bar x,i}^{(\sigma)})^2\right)h-
$$
$$
-\frac{2}{h^2}\sum\limits_{i=1}^{N-1}(v_{i+1}-v_{i})\left(a_{i+1}y_{\bar
x,i+1}^{(\sigma)}\right)^2h=\frac{4}{h^2}\sum\limits_{i=1}^{N}v_{i}\left(a_{i}y_{\bar
x,i}^{(\sigma)}\right)^2h-\frac{2}{h}\sum\limits_{i=2}^{N}v_{\bar
x,i}\left(a_{i}y_{\bar x,i}^{(\sigma)}\right)^2h-
$$
$$
-\frac{2}{h}\left(v_N(a_Ny_{\bar
x,N}^{(\sigma)})^2+v_1(a_1y_{x,0}^{(\sigma)})^2\right)\leq\frac{4c_2}{h^2}\|\sqrt{v}\sqrt{a}y_{\bar
x}^{(\sigma)}]|_0^2-\frac{2}{h}\sum\limits_{i=1}^{N}v_{\bar
x,i}\left(a_{i}y_{\bar x,i}^{(\sigma)}\right)^2h-
$$
$$
-\frac{2}{h}\left(v_N(a_Ny_{\bar
x,N}^{(\sigma)})^2+v_0(a_1y_{x,0}^{(\sigma)})^2\right).
$$

This implies the desired inequality (\ref{ur32}).

{\bf Lemma 4.} The inequality
$$
\|\sqrt{v}\sqrt{a}y_{\bar
x}^{(\sigma)}]|_0^2\geq\frac{h^2}{4c_2(1+\varepsilon)}\|\sqrt{v}\Delta_{0t}^{\gamma}
y\|_0^2+ \frac{h}{2c_2}\sum\limits_{i=1}^{N}v_{\bar x,i}(a_iy_{\bar
x,i}^{(\sigma)})^2h+
$$
\begin{equation}
+\frac{h}{2c_2}\left(v_N(a_Ny_{\bar
x,N}^{(\sigma)})^2+v_0(a_1y_{x,0}^{(\sigma)})^2\right)-\frac{h^2}{4c_2\varepsilon}\|\sqrt{v}\varphi\|_0^2,
\quad \varepsilon>0.
 \label{ur32.5}
\end{equation}
holds for any nonnegative function $v(x)$ defined on the grid $\bar
\omega_h$ and any solution $y(x,t)$ of equation (\ref{ur4}).

{\bf Proof.} Since $y(x,t)$ is a solution of equation (\ref{ur4}),
it follows that $\sqrt{v}\Delta_{0t}^{\gamma}
y_i=\sqrt{v_i}(ay_{\bar x}^{(\sigma)})_{x,i}+\sqrt{v_i}\varphi_i$
for all $i=1,2,\ldots,N-1$ for any nonnegative function $v(x)$;
consequently,

$$
\|\sqrt{v}\Delta_{0t}^{\gamma} y\|_0^2=\|\sqrt{v}(ay_{\bar
x}^{(\sigma)})_{x}+\sqrt{v}\varphi\|_0^2\leq
\frac{1+\varepsilon}{h^2}\sum\limits_{i=1}^{N-1}v_i\left(a_{i+1}y_{\bar
x,i+1}^{(\sigma)}-a_{i}y_{\bar x,i}^{(\sigma)}\right)^2h+
$$
$$
+\left(1+\frac{1}{\varepsilon}\right)\sum\limits_{i=1}^{N-1}v_i\varphi_i^2h\leq\frac{2(1+\varepsilon)}{h^2}\sum\limits_{i=1}^{N-1}v_i\left((a_{i+1}y_{\bar
x,i+1}^{(\sigma)})^2+(a_{i}y_{\bar
x,i}^{(\sigma)})^2\right)h+\left(1+\frac{1}{\varepsilon}\right)\|\sqrt{v}\varphi\|_0^2
=
$$
$$
=\frac{4c_2(1+\varepsilon)}{h^2}\|\sqrt{v_i}\sqrt{a}y_{\bar
x,i}^{(\sigma)}]|_0^2-\frac{2(1+\varepsilon)}{h}\sum\limits_{i=1}^{N}v_{\bar
x,i}\left(a_{i}y_{\bar x,i}^{(\sigma)}\right)^2h-
$$

$$
-\frac{2(1+\varepsilon)}{h}\left(v_N(a_Ny_{\bar
x,N}^{(\sigma)})^2+v_0(a_1y_{x,0}^{(\sigma)})^2\right)+\frac{1+\varepsilon}{\varepsilon}\|\sqrt{v_i}\varphi_i\|_0^2.
$$
Hence we derive the desired inequality (\ref{ur32.5}).

{\bf Theorem 4.} If
$$
(\beta\alpha^{-1}-1)(\alpha^2-1)\geq 0,
$$ then the condition
\begin{equation}
\sigma\geq\frac{1}{3-2^{1-\gamma}}-\frac{h^2(2-2^{1-\gamma})}{2c_2\tau^\gamma(3-2^{1-\gamma})\Gamma(2-\gamma)}
 \label{ur32.7777}
\end{equation} is sufficient
for the stability of the difference scheme (\ref{ur4})--(\ref{ur6})
for $\varphi(x,t)\equiv 0$ and $\mu(t)\equiv 0$, its solution
satisfies the estimate
\begin{equation}
\|y^{n+1}\|_1\leq\|y^0\|_1,
 \label{ur33}
\end{equation}
here
$$\|y\|_1^2=\|y\|_0^2+\delta_1(\alpha,\beta)\|p_1y\|_0^2+\gamma_1(\alpha,\beta)y_0^2h,$$
where
$$
p_1^2(x)=\sum_{s=i}^{N-1}h/a_{s+1}, \quad
\delta_1(\alpha,\beta)=(\beta\alpha^{-1}-1)/p_1^2(0),\quad
\gamma_1(\alpha,\beta)=(\alpha\beta+1)/(2\alpha^2),
$$
 if
$\beta\alpha^{-1}-1\geq0$, $\alpha^2-1\geq0$;\quad
$\|y\|_1^2=\|y\|_0^2+\delta_1(\alpha^{-1},\beta^{-1})\|p_1(1-x)y\|_0^2+\gamma_1(\alpha^{-1},\beta^{-1})\alpha^{-2}y_0^2h$,
if  $\beta\alpha^{-1}-1\leq0$, $\alpha^2-1\leq0$. The norm
$\|y\|_1^2$ is equivalent to the norm
$|[y]|_0^2=0.5hy_0^2+0.5hy_N^2+\|y\|_0^2$.

 {\bf Proof.} By multiplying equation (\ref{ur4}) by $y^{(\sigma)}h$ and by summing with respect to
$i$ from $1$ to $N-1$, we obtain the relation
\begin{equation}
(\Delta_{0t}^{\gamma} y,y^{(\sigma)})-((ay_{\bar
x}^{(\sigma)})_x,y^{(\sigma)})=(\varphi,y^{(\sigma)}),
 \label{ur34}
\end{equation}
where $(y,v)=\sum_{i=1}^{N-1}y_iv_ih$.

Let us transform the terms in relation (\ref{ur34}),

$$ (\Delta_{0t}^{\gamma} y,y^{(\sigma)})=\frac{1}{2}\Delta_{0t}^{\gamma}\left(\|y\|_0^2\right)+
\frac{\tau^\gamma\Gamma(2-\gamma)}{2(2-2^{1-\gamma})}\left((3-2^{1-\gamma})\sigma-1\right)\|\Delta_{0t}^{\gamma}
y\|_0^2+\|\sqrt{J^{(\sigma)}(y)}\|_0^2,
$$
$$
-((ay_{\bar x}^{(\sigma)})_x,y^{(\sigma)})=(ay_{\bar
x}^{(\sigma)},y_{\bar x}^{(\sigma)}]-a_Ny_{\bar
x,N}^{(\sigma)}y_{N}^{(\sigma)}+a_1y_{
x,0}^{(\sigma)}y_{0}^{(\sigma)}=
$$
$$
=\|\sqrt{a}y_{\bar
x}^{(\sigma)}]|_0^2+\left(1-\frac{\beta}{\alpha}\right)a_1y_{
x,0}^{(\sigma)}y_{0}^{(\sigma)}-
y_N^{(\sigma)}\mu(t_{n+\sigma})-\frac{h}{2}(\varphi_N+\beta\varphi_0)y_N^{(\sigma)}+
$$
\begin{equation}
+
\frac{\alpha\beta+1}{2\alpha^2}\left(\frac{1}{2}\Delta_{0t}^{\gamma}(y_0^2)
+\frac{\tau^\gamma\Gamma(2-\gamma)}{2(2-2^{1-\gamma})}\left((3-2^{1-\gamma})\sigma-1\right)(\Delta_{0t}^{\gamma}y_{t,0})^2+J^{(\sigma)}(y_0)\right)h,
 \label{ur35}
\end{equation}
where $J^{(\sigma)}(y)=\sigma J_1(y)+(1-\sigma)J_2(y)$.

By substituting the expression (\ref{ur35}) into relation
(\ref{ur34}), we obtain

$$
\frac{1}{2}\Delta_{0t}^{\gamma}\left(\|y\|_2^2\right)+
\frac{\tau^\gamma\Gamma(2-\gamma)}{2(2-2^{1-\gamma})}\left((3-2^{1-\gamma})\sigma-1\right)\|\Delta_{0t}^{\gamma}
y\|_2^2+\|\sqrt{a}y_{\bar
x}^{(\sigma)}]|_0^2+\|\sqrt{J^{(\sigma)}(y)}\|_2^2 =
$$
\begin{equation}
=\left(\frac{\beta}{\alpha}-1\right)a_1y_{
x,0}^{(\sigma)}y_{0}^{(\sigma)}+y_N^{(\sigma)}\mu(t_{n+\sigma})+\frac{h}{2}(\varphi_N+\beta\varphi_0)y_N^{(\sigma)}+(\varphi,y^{(\sigma)}),
 \label{ur36}
\end{equation}
where $\|y\|_2^2=\|y\|_0^2+(\alpha\beta+1)/(2\alpha^2)y_0^2h$.

By multiplying equation (\ref{ur4}) by $y_{i}^{(\sigma)}h$ and by
summing with respect to $s$ from $1$ to $i$, we obtain

\begin{equation}
\sum\limits_{s=1}^{i}y_{s}^{(\sigma)}\Delta_{0t}^{\gamma}y_{s}h-\sum\limits_{s=1}^{i}(ay_{\bar
x}^{(\sigma)})_{x,s}y_{s}^{(\sigma)}h=\sum\limits_{s=1}^{i}\varphi_sy_{s}^{(\sigma)}h,\quad
i=0,1,...,N-1.
 \label{ur37}
\end{equation}
Here we adopt the convention that the sum is zero if the upper limit
of the summation is less than the lower one.

Let us transform the terms in relation (\ref{ur37}),
$$
\sum\limits_{s=1}^{i}y_{s}^{(\sigma)}\Delta_{0t}^{\gamma}y_{s}h=
\frac{1}{2}\Delta_{0t}^{\gamma}\left(\sum\limits_{s=1}^{i}y_{s}^2h\right)+
$$
$$
+\frac{\tau^\gamma\Gamma(2-\gamma)}{2(2-2^{1-\gamma})}\left((3-2^{1-\gamma})\sigma-1\right)
\sum\limits_{s=1}^{i}\left(\Delta_{0t}^{\gamma}y_{s}\right)^2h+\sum\limits_{s=1}^{i}J^{(\sigma)}(y_i)h,
$$

$$
-\sum\limits_{s=1}^{i}(ay_{\bar
x}^{(\sigma)})_{x,s}y_{s}^{(\sigma)}h=
\sum\limits_{s=1}^{i+1}a_s(y_{\bar
x,s}^{(\sigma)})^2h-a_{i+1}y_{\bar
x,i+1}^{(\sigma)}y_{i+1}^{(\sigma)}+a_{1}y_{0
x,0}^{(\sigma)}y_{0}^{(\sigma)}=
$$

\begin{equation}
=\sum\limits_{s=1}^{i}a_s(y_{\bar
x,s}^{(\sigma)})^2h+\frac{h}{2}a_{i+1}(y_{\bar x,i+1}^{(\sigma)})^2-
\frac{1}{2}a_{i+1}\left((y^{(\sigma)})^2\right)_{\bar
x,i+1}+a_{1}y_{0 x,0}^{(\sigma)}y_{0}^{(\sigma)}. \label{ur38}
\end{equation}

By substituting the expressions (\ref{ur38}) into relation
(\ref{ur37}), we get

$$
\frac{1}{2}\Delta_{0t}^{\gamma}\left(\sum\limits_{s=1}^{i}y_{s}^2h\right)
+\frac{\tau^\gamma\Gamma(2-\gamma)}{2(2-2^{1-\gamma})}\left((3-2^{1-\gamma})\sigma-1\right)
\sum\limits_{s=1}^{i}\left(\Delta_{0t}^{\gamma}y_{s}\right)^2h+
$$
$$
+\sum\limits_{s=1}^{i}a_s(y_{\bar
x,s}^{(\sigma)})^2h+\frac{h}{2}a_{i+1}(y_{\bar
x,i+1}^{(\sigma)})^2+\sum\limits_{s=1}^{i}J^{(\sigma)}(y_i)h=
$$
\begin{equation}
=\frac{1}{2}a_{i+1}\left((y^{(\sigma)})^2\right)_{\bar
x,i+1}+a_{1}y_{0
x,0}^{(\sigma)}y_{0}^{(\sigma)}+\sum\limits_{s=1}^{i}\varphi_sy_{s}^{(\sigma)}h.
\label{ur39}
\end{equation}

If we multiply relation (\ref{ur39}) by $h/a_{i+1}$  and sum with
respect to $i$ from $0$ to $N-1$, then we obtain

$$
\frac{1}{2}\Delta_{0t}^{\gamma}\left(\|p_1y\|_0^2\right)+\frac{\tau^\gamma\Gamma(2-\gamma)}{2(2-2^{1-\gamma})}\left((3-2^{1-\gamma})\sigma-1\right)
\|p_1\Delta_{0t}^{\gamma}y\|_0^2+
$$
$$
+\|p_1\sqrt{a}(y_{\bar x}^{(\sigma)})\|_0^2+\frac{h}{2}\|y_{\bar
x}^{(\sigma)}]|_0^2+\|p_1\sqrt{J^{(\sigma)}(y)}\|_0^2=
$$
\begin{equation}
=\frac{1}{2}(1-\alpha^2)(y_N^{(\sigma)})^2+p_1^2(0)a_{1}y_{
x,0}^{(\sigma)}y_{0}^{(\sigma)}+(p_1^2\varphi,y^{(\sigma)}),
\label{ur40}
\end{equation}
where $p_1(x_i)=\left(\sum_{s=i}^{N-1}h/a_{s+1}\right)^{1/2}$,
$p_1(x_N)=0$.

Let us multiply relation (\ref{ur40}) by
$\delta_1=(\beta\alpha^{-1}-1)/p_1^2(0)$, and add the resulting
relation to identity (\ref{ur36}),
$$
\frac{1}{2}\Delta_{0t}^{\gamma}\left(\|y\|_2^2+\delta_1\|p_1y\|_0^2\right)+\frac{\tau^\gamma\Gamma(2-\gamma)}{2(2-2^{1-\gamma})}\left((3-2^{1-\gamma})\sigma-1\right)
(\|\Delta_{0t}^{\gamma}y\|_2^2+\delta_1\|p_1\Delta_{0t}^{\gamma}y\|_0^2)+
$$
$$
+\|\sqrt{a}y_{\bar x}^{(\sigma)}]|_0^2+\delta_1\|p_1\sqrt{a}y_{\bar
x}^{(\sigma)}\|_0^2+\frac{\delta_1h}{2}\|y_{\bar x}^{(\sigma)}]|_0^2
+\frac{\delta_1}{2}(\alpha^2-1)(y_N^{(\sigma)})^2+\|\sqrt{J^{(\sigma)}(y)}\|_2^2+
$$
\begin{equation}
+\delta_1\|p_1\sqrt{J^{(\sigma)}(y)}\|_0^2=((1+\delta_1p_1^2)\varphi,y^{(\sigma)})+\frac{h}{2}(\varphi_N+\beta\varphi_0)y_N^{(\sigma)}+y_N^{(\sigma)}\mu(t_{n+\sigma}).
\label{ur41}
\end{equation}

Consider two cases: {\bf 1)} $\beta\alpha^{-1}-1\geq0$ and
$\alpha^2-1\geq0$;  {\bf 2)} $\beta\alpha^{-1}-1\leq0$ and
$\alpha^2-1\leq0$.

{\bf 1)} Let $\varphi$ and $\mu\equiv0$. For $v(x)=1$, $x\in \bar
\omega_h$ and $v(x)=p_1^2(x)$, $i=0,1,\ldots,N$ from Lemma 3 we
obtain the inequalities

\begin{equation}
\|\sqrt{a}y_{\bar
x}^{(\sigma)}]|_0^2\geq\frac{h^2}{4c_2}\|\Delta_{0t}^{\gamma}y\|_0^2
+\frac{h}{2c_2}\left((a_Ny_{\bar
x,N}^{(\sigma)})^2+(a_1y_{x,0}^{(\sigma)})^2\right),
 \label{ur42}
\end{equation}

\begin{equation}
\|p_1\sqrt{a}y_{\bar
x}^{(\sigma)}]|_0^2\geq\frac{h^2}{4c_2}\|p_1\Delta_{0t}^{\gamma}y\|_0^2
-\frac{h}{2c_2}\sum\limits_{i=1}^{N}a_i(y_{\bar
x,i}^{(\sigma)})^2h+\frac{h}{2c_2}p_1^2(0)(a_1y_{x,0}^{(\sigma)})^2.
 \label{ur43}
\end{equation}

Relation (\ref{ur41}), together with (\ref{ur42}) and (\ref{ur43}),
implies the inequality

$$
\frac{1}{2}\Delta_{0t}^{\gamma}\left(\|y\|_2^2+\delta_1\|p_1y\|_0^2\right)
+\left(\frac{\tau^\gamma\Gamma(2-\gamma)}{2(2-2^{1-\gamma})}\left((3-2^{1-\gamma})\sigma-1\right)+\frac{h^2}{4c_2}\right)
\|\Delta_{0t}^{\gamma}y\|_0^2+
$$
$$
+\delta_1\left(\frac{\tau^\gamma\Gamma(2-\gamma)}{2(2-2^{1-\gamma})}\left((3-2^{1-\gamma})\sigma-1\right)+\frac{h^2}{4c_2}\right)\|p_1y_t\|_0^2+
\frac{\delta_1h}{2}\left(\|y_{\bar
x}^{(\sigma)}]|_0^2-\frac{1}{c_2}\sum\limits_{i=1}^{N}a_i(y_{\bar
x,i}^{(\sigma)})^2h\right)+
$$
$$
+\frac{\tau^\gamma\Gamma(2-\gamma)}{2(2-2^{1-\gamma})}\left((3-2^{1-\gamma})\sigma-1\right)\frac{\alpha\beta+1}{2\alpha^2}\left(\Delta_{0t}^{\gamma}y_{0}\right)^2h+
\frac{h}{2c_2}\left((a_Ny_{\bar
x,N}^{(\sigma)})^2+(a_1y_{x,0}^{(\sigma)})^2\right)+
$$
\begin{equation}
+\frac{h}{2c_2}\delta_1p_1^2(0)(a_1y_{x,0}^{(\sigma)})^2
+\frac{\delta_1}{2}(\alpha^2-1)(y_N^{(\sigma)})^2+\|\sqrt{J^{(\sigma)}(y)}\|_2^2+\delta_1\|p_1\sqrt{J^{(\sigma)}(y)}\|_0^2\leq
0. \label{ur44}
\end{equation}

Since
$$\sigma\geq\frac{1}{3-2^{1-\gamma}}-\frac{h^2(2-2^{1-\gamma})}{2c_2\tau^\gamma(3-2^{1-\gamma})\Gamma(2-\gamma)},\quad
\delta_1\left(\|y_{\bar
x,i+1}^{(\sigma)}]|_0^2-\frac{1}{c_2}\sum\limits_{i=1}^{N}a_i(y_{\bar
x,i}^{(\sigma)})^2h\right)\geq 0,
$$
$$
\delta_1\geq0,\quad \frac{\delta_1}{2}(\alpha^2-1)\geq 0, \quad
\|\sqrt{J^{(\sigma)}(y)}\|_2^2+\delta_1\|p_1\sqrt{J^{(\sigma)}(y)}\|_0^2\geq
0,
$$ it follows from (\ref{ur44}) that
$$
\Delta_{0t}^{\gamma}\left(\|y\|_2^2+\delta_1\|p_1y\|_0^2\right)
+\frac{\tau^\gamma\Gamma(2-\gamma)}{2(2-2^{1-\gamma})}\left((3-2^{1-\gamma})\sigma-1\right)\frac{\alpha\beta+1}{\alpha^2}\left(\Delta_{0t}^{\gamma}y_{0}\right)^2h+
$$
\begin{equation}
+\frac{h}{c_2}\left((a_Ny_{\bar
x,N}^{(\sigma)})^2+\frac{\beta}{\alpha}(a_1y_{x,0}^{(\sigma)})^2\right)
\leq 0. \label{ur45}
\end{equation}

If $\sigma\geq 1/(3-2^{1-\gamma})$, $\varphi\equiv0$ and
$\mu\equiv0$, then from (\ref{ur45}) follows inequality
\begin{equation}
\Delta_{0t}^{\gamma}\|y\|_1^2\leq0. \label{ur45.55}
\end{equation}
where $\|y\|_1=\left(\|y\|_2^2+\delta_1\|p_1y\|_0^2\right)^{1/2}.$

Consider the case in which $0\leq \sigma<1/(3-2^{1-\gamma})$. We
introduce the notation
$$
\xi=\frac{\tau^\gamma\Gamma(2-\gamma)}{2(2-2^{1-\gamma})}\left((3-2^{1-\gamma})\sigma-1\right).
$$
By virtue of the condition (\ref{ur32.7777}), in this case, we have
$-h^2/(4c_2)\leq \xi<0$.

It follows from the boundary conditions (\ref{ur5}) with $\varphi$
and $\mu\equiv0$ that
$$
\frac{\alpha\beta+1}{\alpha}\Delta_{0t}^{\gamma}y_{0}=-\frac{2}{h}\left(a_Ny_{\bar
x, N}^{(\sigma)}-\beta a_1 y_{ x, 0}^{(\sigma)}\right).
$$

We square both sides of the last relation and divide the resulting
relation by $(\alpha\beta+1)$,
$$
\frac{\alpha\beta+1}{\alpha^2}\left(\Delta_{0t}^{\gamma}y_{0}\right)^2=
\frac{4}{h^2(\alpha\beta+1)}(a_Ny_{\bar x,N}^{(\sigma)})^2+
$$
\begin{equation}
+\frac{8\beta}{h^2(\alpha\beta+1)}(a_Ny_{\bar
x,N}^{(\sigma)})(a_1y_{x,0}^{(\sigma)})+\frac{4\beta^2}{h^2(\alpha\beta+1)}
(a_1y_{x,0}^{(\sigma)})^2. \label{ur46}
\end{equation}

Inequality (\ref{ur45}), together with (\ref{ur46}), implies that
$$
\Delta_{0t}^{\gamma}\left(\|y\|_2^2+\delta_1\|p_1y\|_0^2\right)+\left(\frac{4\xi}{(\alpha\beta+1)
h}+\frac{h}{c_2}\right)(a_Ny_{\bar x,N}^{(\sigma)})^2-
$$
\begin{equation}
-\frac{8\beta\xi}{(\alpha\beta+1) h}(a_Ny_{\bar
x,N}^{(\sigma)})(a_1y_{x,0}^{(\sigma)})+
\left(\frac{4\beta^2\xi}{(\alpha\beta+1) h} +\frac{h\beta}{\alpha
c_2}\right)(a_1y_{x,0}^{(\sigma)})^2\leq 0.
 \label{ur47}
\end{equation}

Note that
 $$
\frac{4\xi}{(\alpha\beta+1) h}+\frac{h}{c_2}\geq
\frac{h}{c_2}-\frac{h}{(\alpha\beta+1)
c_2}=\frac{h\alpha\beta}{(\alpha\beta+1) c_2}>0,
 $$

The quadratic form
$$
\left(\frac{4\xi}{(\alpha\beta+1) h}+\frac{h}{c_2}\right)(a_Ny_{\bar
x,N}^{(\sigma)})^2-\frac{8\beta\xi}{(\alpha\beta+1) h}(a_Ny_{\bar
x,N}^{(\sigma)})(a_1y_{x,0}^{(\sigma)})+
\left(\frac{4\beta^2\xi}{(\alpha\beta+1) h} +\frac{h\beta}{\alpha
c_2}\right)(a_1y_{x,0}^{(\sigma)})^2
$$
is nonnegative for all values of $a_Ny_{\bar x,N}^{(\sigma)}$ and
$a_0y_{x,0}^{(\sigma)}$ if and only if
$$
\frac{16\beta^2\xi^2}{(\alpha\beta+1)^2h^2}-
\left(\frac{4\xi}{(\alpha\beta+1) h}+\frac{h}{c_2}\right)
\left(\frac{4\beta^2\xi}{(\alpha\beta+1) h} +\frac{\beta h}{\alpha
c_2}\right)\leq 0,
$$
which, after simple transformations, acquires the form
$$
-\frac{4\beta}{\alpha c_2
}\left(\frac{\tau^\gamma\Gamma(2-\gamma)}{2(2-2^{1-\gamma})}\left((3-2^{1-\gamma})\sigma-1\right)+\frac{h^2}{4c_2}\right)\leq
0.
$$
Consequently, the estimate (\ref{ur45.55}) follows from inequality
(\ref{ur47}).

Let us rewrite (\ref{ur45.55}) as
\begin{equation}
\frac{1}{\Gamma(2-\gamma)}\sum\limits_{s=0}^{n}(t_{n-s+1}^{1-\alpha}-t_{n-s}^{1-\alpha})\frac{\|y^{s+1}\|_1^2-\|y^{s}\|_1^2}{\tau}\leq
0.
 \label{ur47.55}
\end{equation}
It is obvious that at $n=0$ the a priori estimate (\ref{ur33})
follows from (\ref{ur47.55}). Let us prove that (\ref{ur33}) holds
for $n=1,2,\ldots$ by using the mathematical induction method. For
this purpose, let us assume that the a priori estimate (\ref{ur33})
takes place for all $n=0,1,\ldots,k-1$, $k=1,2,\ldots$. From
(\ref{ur47.55}) at $n=k$ one has
\begin{equation}
\tau^{1-\gamma}\|y^{k+1}\|_1^2\leq\sum\limits_{s=1}^{k}(-t_{k-s+2}^{1-\gamma}+2t_{k-s+1}^{1-\gamma}-t_{k-s}^{1-\gamma})\|y^{s}\|_1^2+
(t_{k+1}^{1-\gamma}-t_k^{1-\gamma})\|y^0\|_1^2.
 \label{ur47.66}
\end{equation}
Since $-t_{k+1}^{1-\gamma}+2t_{k}^{1-\gamma}-t_{k-1}^{1-\gamma}>0$
for all $k=1,2,\ldots$ \cite{ShkhTau:06}, and by the assumption of
the mathematical induction $\|y^{s}\|_0^2\leq\|y^0\|_0^2$ at
$s=1,2,\ldots,k$, then from (\ref{ur47.66}) one obtains the
following inequality:
$$
\tau^{1-\gamma}\|y^{k+1}\|_1^2\leq\left(
(t_{k+1}^{1-\gamma}-t_k^{1-\gamma})+\sum\limits_{s=1}^{k}(-t_{k-s+2}^{1-\gamma}+2t_{k-s+1}^{1-\gamma}-t_{k-s}^{1-\gamma})\right)\|y^0\|_1^2=
$$
\begin{equation}
=\left((t_{k+1}^{1-\gamma}-t_k^{1-\gamma})-t_{k+1}^{1-\gamma}+t_k^{1-\gamma}+t_{1}^{1-\gamma}-t_0^{1-\gamma}\right)\|y^0\|_1^2=
\tau^{1-\gamma}\|y^{0}\|_1^2.
 \label{ur47.77}
\end{equation}

{\bf 2)} The a priori estimate (\ref{ur33}) for the second case, in
which $\beta\alpha^{-1}-1\leq0, \alpha^2-1\leq0$, follows directly
from the first case. Indeed, if we set $y(x,t)=v(1-x,t)$, then the
function $v(x,t)$ satisfies the problem
\begin{equation}
\Delta_{0t_n}^{\gamma}v_{i}-(\bar av_{\bar
x}^{(\sigma)})_{x,i}=\bar\varphi_i^n,\quad i=1,2,\ldots,N-1,
\label{ur48}
\end{equation}

\begin{equation}
\begin{cases}
 v_0^{n+1}-\dfrac{1}{\alpha} v_N^{n+1}=0,\\
\dfrac{1}{\beta}\Delta_{0t_n}^{\gamma}
v_{0}+\Delta_{0t_n}^{\gamma}v_{N}+\dfrac{2}{h}\left(\bar a_Nv_{\bar
x, N}^{(\sigma)}-\dfrac{1}{\beta} \bar a_1 v_{ x,
0}^{(\sigma)}\right)=\dfrac{2}{\beta h}\mu(t_{n+1/2})+\bar\varphi_N+\dfrac{1}{\beta}\bar \varphi_0,\\
\end{cases}
\label{ur49}
\end{equation}

\begin{equation}
v_i^0=u_0(1-x_i),
 \label{ur50}
\end{equation}
where $\bar a_i=a_{N-i+1}$, $\bar \varphi_i^n=\varphi_{N-i}^n$.

By virtue of the conditions
$$0<c_1\leq \bar a\leq c_2,\, \alpha\beta^{-1}-1\geq0,\,
\alpha^{-2}-1\geq0,\,
\sigma\geq\frac{1}{3-2^{1-\gamma}}-\frac{h^2(2-2^{1-\gamma})}{2c_2\tau^\gamma(3-2^{1-\gamma})\Gamma(2-\gamma)},
$$
the estimate (\ref{ur33}) holds for the solution $v(x,t)$ of problem
(\ref{ur48})--(\ref{ur50}) for $\varphi\equiv0$ and $\mu\equiv0$.
Consequently, by virtue of the relations $\|v\|_0^2=\|y\|_0^2$,
$\|p_1v\|_0^2=\|p_1(1-x)y\|_0^2$ and $v_0^2=\alpha^{-2}y_0^2$, the
solution of problem (\ref{ur4})--(\ref{ur6}) for $\varphi\equiv0$
and $\mu\equiv0$, satisfies the a priori estimate (\ref{ur33}). The
proof of Theorem 4 is complete.

{\bf Theorem 5.} If
$$
(\beta\alpha^{-1}-1)(\alpha^2-1)> 0,
$$ then the condition
\begin{equation}
\sigma\geq\frac{1}{3-2^{1-\alpha}}-\frac{h^2(2-2^{1-\alpha})(1-\varepsilon)}{2c_2\tau^\alpha(3-2^{1-\alpha})\Gamma(2-\alpha)},
\quad 0<\varepsilon<1
 \label{ur50.55}
\end{equation}
 is sufficient for the stability of the difference scheme
(\ref{ur4})--(\ref{ur6}) and its solution satisfies the estimate
\begin{equation}
\frac{1}{\Gamma(2-\gamma)}\sum\limits_{s=0}^{n}(t_{n-s+1}^{1-\gamma}-t_{n-s}^{1-\gamma})\|y^{s+1}\|_0^2\leq
\frac{t_{n+1}^{1-\gamma}}{\Gamma(2-\gamma)}\|y^{0}\|_0^2+
{M}\left(\sum\limits_{s=0}^{n}(|[\varphi^s]|_0^2+\mu^2(t_{s+\sigma}))\tau\right),
 \label{ur51}
\end{equation}
where ${M}>0$ is a known number independent of $h$, $\tau$ and
$t_n$.

 {\bf Proof.} The solution of the difference scheme
(\ref{ur4})--(\ref{ur6}) satisfies identity (\ref{ur41}).

Consider two cases: {\bf 1)}  $\beta\alpha^{-1}-1>0$, $\alpha^2-1>0$
and {\bf 2)} $\beta\alpha^{-1}-1<0$, $\alpha^2-1<0$.

{\bf 1)} If $v(x)=1$ and $v(x)=p_1^2(x), i=0,1,...,N$, then Lemma 4
readily implies the inequalities

\begin{equation}
\|\sqrt{a}y_{\bar
x}^{(\sigma)}]|_0^2\geq\frac{h^2}{4c_2(1+\varepsilon_1)}\|\Delta_{0t}^{\gamma}y\|_0^2
+\frac{h}{2c_2}\left((a_Ny_{\bar
x,N}^{(\sigma)})^2+(a_1y_{x,0}^{(\sigma)})^2\right)-\frac{h^2}{4c_2\varepsilon_1}\|\varphi\|_0^2,
 \label{ur52}
\end{equation}

$$
\|p_1\sqrt{a}y_{\bar
x}^{(\sigma)}]|_0^2\geq\frac{h^2}{4c_2(1+\varepsilon_1)}\|p_1\Delta_{0t}^{\gamma}y\|_0^2-
$$
\begin{equation}
-\frac{h}{2c_2}\sum\limits_{i=1}^{N}a_i(y_{\bar
x,i}^{(\sigma)})^2h+\frac{h}{2c_2}p_1^2(0)(a_1y_{x,0}^{(\sigma)})^2-\frac{h^2}{4c_2\varepsilon_1}\|p_1\varphi\|_0^2.
 \label{ur53}
\end{equation}

Let us estimate the right-hand side of identity (\ref{ur41}). By
virtue of the relation
$y_i^{(\sigma)}=y_N^{(\sigma)}-\sum\limits_{s=i+1}^{N}y_{\bar x,
s}^{(\sigma)}h$, $i=0,1,...,N-1$, we have
$$
((1+\delta_1p_1^2)\varphi,y^{(\sigma)})=
y_N^{(\sigma)}\sum\limits_{i=1}^{N-1}(1+\delta_1p_1^2(x_i))\varphi_ih-
\sum\limits_{i=1}^{N-1}(1+\delta_1p_1^2(x_i))\varphi_ih\sum\limits_{s=i+1}^{N}y_{\bar
x,s}^{(\sigma)}h\leq
$$
$$
\leq\frac{\varepsilon_2}{2}(y_N^{(\sigma)})^2+\frac{\bar
\gamma_1}{2\varepsilon_2}\sum\limits_{i=1}^{N-1}(1+\delta_1p_1^2(x_i))\varphi_i^2h-
$$
$$-
\sum\limits_{i=1}^{N-1}\sqrt{(1+\delta_1p_1(x_{i+1}))a_{i+1}}y_{x,i}^{(\sigma)}h
\sum\limits_{s=1}^{i}\frac{(1+\delta_1p_1^2(x_s))}{\sqrt{(1+\delta_1p_1(x_{i+1}))a_{i+1}}}\varphi_sh\leq
$$
\begin{equation}
\leq\frac{\varepsilon_2}{2}(y_N^{(\sigma)})^2+\frac{\varepsilon_1}{2}\sum\limits_{i=2}^{N}(1+\delta_1p_1(x_{i}))a_{i}(y_{\bar
x,i}^{(\sigma)})^2h + \left(\frac{\bar
\gamma_1}{2\varepsilon_2}+\frac{\bar \gamma_1\bar
\gamma_2}{2\varepsilon_1}\right)\sum\limits_{i=1}^{N-1}(1+\delta_1p_1^2(x_i))\varphi_i^2h,
 \label{ur54}
\end{equation}
where $\varepsilon_1$,  $\varepsilon_2>0$, $\bar
\gamma_1=\sum_{i=1}^{N-1}(1+\delta_1p_1^2(x_i))h$, $\bar
\gamma_2=\sum_{i=1}^{N-1}(1+\delta_1p_1^2(x_i))^{-1}a_i^{-1}h$, and
\begin{equation}
\frac{h}{2}(\varphi_N+\beta\varphi_0)y_N^{(\sigma)}+y_N^{(\sigma)}\mu(t_{n+\sigma})
\leq\frac{\varepsilon_2}{2}(y_N^{(\sigma)})^2+\frac{1}{2\varepsilon_2}\bar\mu^2(t_{n+\sigma}),
 \label{ur55}
\end{equation}
where
$\bar\mu(t_{n+\sigma})=\mu(t_{n+\sigma})+0.5h(\varphi_N+\beta\varphi_0)$.

Relation (\ref{ur41}), together with inequalities (\ref{ur54}) and
(\ref{ur55}), implies the inequality

$$
\frac{1}{2}\Delta_{0t}^{\gamma}\left(\|y\|_2^2+\delta_1\|p_1y\|_0^2\right)+\frac{\tau^\gamma\Gamma(2-\gamma)}{2(2-2^{1-\gamma})}\left((3-2^{1-\gamma})\sigma-1\right)
(\|\Delta_{0t}^{\gamma}y\|_2^2+\delta_1\|p_1\Delta_{0t}^{\gamma}y\|_0^2)+
$$
$$
+(1-\varepsilon_1)\|\sqrt{a}y_{\bar
x}^{(\sigma)}]|_0^2+\delta_1(1-\varepsilon_1)\|p_1\sqrt{a}y_{\bar
x}^{(\sigma)}\|_0^2+\frac{\delta_1h}{2}\|y_{\bar x}^{(\sigma)}]|_0^2
+\left(\frac{\delta_1}{2}(\alpha^2-1)-\varepsilon_2\right)(y_N^{(\sigma)})^2+
$$
\begin{equation}
+\|\sqrt{J^{(\sigma)}(y)}\|_2^2+\delta_1\|p_1\sqrt{J^{(\sigma)}(y)}\|_0^2\leq
M_2(\varepsilon_1,\varepsilon_2)\left(|[\varphi]|_0^2+\mu^2(t_{n+\sigma})\right).
\label{ur56}
\end{equation}

By taking into account inequalities (\ref{ur52}) and (\ref{ur53}),
from inequality (\ref{ur56}) with
$\varepsilon_2=\delta_1(\alpha^2-1)/2$, we obtain the inequality

$$
\frac{1}{2}\Delta_{0t}^{\gamma}\left(\|y\|_2^2+\delta_1\|p_1y\|_0^2\right)
+\left(\frac{\tau^\gamma\Gamma(2-\gamma)}{2(2-2^{1-\gamma})}\left((3-2^{1-\gamma})\sigma-1\right)+\frac{h^2(1-\varepsilon_1)}{4c_2(1+\varepsilon_1)}\right)
\|\Delta_{0t}^{\gamma}y\|_0^2+
$$
$$
+\delta_1\left(\frac{\tau^\gamma\Gamma(2-\gamma)}{2(2-2^{1-\gamma})}\left((3-2^{1-\gamma})\sigma-1\right)+\frac{h^2(1-\varepsilon_1)}{4c_2(1+\varepsilon_1)}\right)\|p_1y_t\|_0^2+
$$
$$
+\frac{\delta_1h}{2}\left(\|y_{\bar
x}^{(\sigma)}]|_0^2-\frac{1}{c_2}\sum\limits_{i=1}^{N}a_i(y_{\bar
x,i}^{(\sigma)})^2h\right)
+\frac{\tau^\gamma\Gamma(2-\gamma)}{2(2-2^{1-\gamma})}\left((3-2^{1-\gamma})\sigma-1\right)\frac{\alpha\beta+1}{2\alpha^2}\left(\Delta_{0t}^{\gamma}y_{0}\right)^2h+
$$
\begin{equation}
+\frac{h}{2c_2}\left((a_Ny_{\bar
x,N}^{(\sigma)})^2+(a_1y_{x,0}^{(\sigma)})^2\right)+\frac{h}{2c_2}\delta_1p_1^2(0)(a_1y_{x,0}^{(\sigma)})^2
\leq
M_3(\varepsilon_1)\left(|[\varphi]|_0^2+\mu^2(t_{n+\sigma})\right).
\label{ur57}
\end{equation}
where $0<\varepsilon_1<1$.

Then we have $\varepsilon=2\varepsilon_1/(1+\varepsilon_1)$,
$0<\varepsilon<1$.

Inequality (\ref{ur57}), together with the assumptions of Theorem 5,
implies the inequality
$$ \left(\|y\|_2^2+\delta_1\|p_1y\|_0^2\right)_t
+\tau\left(\sigma-\frac{1}{2}\right)\frac{\alpha\beta+1}{\alpha^2}y_{t,0}^2h
+
$$
\begin{equation}
+\frac{h}{c_2}\left((a_Ny_{\bar
x,N}^{(\sigma)})^2+\frac{\beta}{\alpha}(a_1y_{x,0}^{(\sigma)})^2\right)
\leq
M_4(\varepsilon)\left(|[\varphi]|_0^2+\mu^2(t_{n+\sigma})\right).
\label{ur58}
\end{equation}

If $\sigma\geq 1/(3-2^{1-\gamma})$, then from (\ref{ur58}) we have
\begin{equation}
\Delta_{0t}^{\gamma}\|y\|_1^2\leq
M\left(|[\varphi]|_0^2+\mu^2(t_{n+\sigma})\right). \label{ur58.55}
\end{equation}
where $\|y\|_1=\left(\|y\|_2^2+\delta_1\|p_1y\|_0^2\right)^{1/2}$,
$M=M_4(\varepsilon)$.

Consider the case in which $0\leq \sigma<1/(3-2^{1-\gamma})$. We
introduce the notation
$$
\xi=\frac{\tau^\gamma\Gamma(2-\gamma)}{2(2-2^{1-\gamma})}\left((3-2^{1-\gamma})\sigma-1\right).
$$
In this case, by virtue of the condition (\ref{ur50.55}), we have
$-h^2(1-\varepsilon)/(4c_2)\leq \xi<0$.

It follows from the boundary conditions (\ref{ur5}) that

$$
\frac{\alpha\beta+1}{\alpha}\Delta_{0t}^{\gamma}y_{0}=-\frac{2}{h}\left(a_Ny_{\bar
x, N}^{(\sigma)}-\beta a_1 y_{ x,
0}^{(\sigma)}\right)+\dfrac{2}{h}\mu(t_{n+\sigma})+\varphi_N+\beta\varphi_0
$$

whence we obtain the inequality
$$
\frac{\alpha\beta+1}{\alpha^2}\Delta_{0t}^{\gamma}y_{0}^2\leq(1+\varepsilon_3)\frac{4}{(\alpha\beta+1)h^2}\left(a_Ny_{\bar
x, N}^{(\sigma)}-\beta a_1 y_{ x,
0}^{(\sigma)}\right)^2+(1+\frac{1}{\varepsilon_3})\frac{4}{h^2}\bar
\mu^2(t_{n+\sigma}), \, \varepsilon_3>0,
$$
or (after the multiplication of the last inequality by $\xi<0$)
\begin{equation}
\frac{\alpha\beta+1}{\alpha^2}\xi
\Delta_{0t}^{\gamma}y_{0}^2\geq\frac{4\xi}{(1-\varepsilon)(\alpha\beta+1)h^2}\left(a_Ny_{\bar
x, N}^{(\sigma)}-\beta a_1 y_{ x,
0}^{(\sigma)}\right)^2-\frac{1-\varepsilon}{c_2\varepsilon}\bar
\mu^2(t_{n+\sigma}), \label{ur59}
\end{equation}
where $\varepsilon=\varepsilon_3/(1+\varepsilon_3)$,
$0<\varepsilon<1$ for $\varepsilon_3>0$.

From inequalities (\ref{ur58}) and (\ref{ur59}),  we have
$$
\Delta_{0t}^{\gamma}\left(\|y\|_2^2+\delta_1\|p_1y\|_0^2\right)+\left(\frac{4\xi}{(1-\varepsilon)(\alpha\beta+1)
h}+\frac{h}{c_2}\right)(a_Ny_{\bar x,N}^{(\sigma)})^2-
$$
$$
 -\frac{8\beta\xi}{(1-\varepsilon)(\alpha\beta+1) h}(a_Ny_{\bar
x,N}^{(\sigma)})(a_1y_{x,0}^{(\sigma)})+
\left(\frac{4\beta^2\xi}{(1-\varepsilon)(\alpha\beta+1) h}
+\frac{h\beta}{\alpha c_2}\right)(a_1y_{x,0}^{(\sigma)})^2\leq
$$
\begin{equation}
\leq M_5(\varepsilon)\left(|[\varphi]|_0^2+\mu^2(t_{n+1/2})\right).
 \label{ur60}
\end{equation}

The quadratic form
$$
\left(\frac{4\xi}{(1-\varepsilon)(\alpha\beta+1)
h}+\frac{h}{c_2}\right)(a_Ny_{\bar x,N}^{(\sigma)})^2-
\frac{8\beta\xi}{(1-\varepsilon)(\alpha\beta+1) h}(a_Ny_{\bar
x,N}^{(\sigma)})(a_1y_{x,0}^{(\sigma)})+
$$
$$
+\left(\frac{4\beta^2\xi}{(1-\varepsilon)(\alpha\beta+1) h}
+\frac{h\beta}{\alpha c_2}\right)(a_1y_{x,0}^{(\sigma)})^2
$$
is nonnegative, because
$$
\frac{4\xi}{(1-\varepsilon)(\alpha\beta+1) h}+\frac{h}{c_2}\geq
\frac{h}{c_2}-\frac{h}{(\alpha\beta+1)c_2}=\frac{\alpha\beta
h}{\alpha\beta+1}>0,
$$
$$
\frac{16\beta^2\xi^2}{(1-\varepsilon)^2(\alpha\beta+1)^2h^2}-
\left(\frac{4\xi}{(1-\varepsilon)(\alpha\beta+1)
h}+\frac{h}{c_2}\right)
\left(\frac{4\beta^2\xi}{(1-\varepsilon)(\alpha\beta+1) h}
+\frac{\beta h}{\alpha c_2}\right)\leq 0.
$$
The last inequality is equivalent to the relation
$$
-\frac{4\beta}{\alpha(1-\varepsilon) c_2
}\left(\frac{\tau^\gamma\Gamma(2-\gamma)}{2(2-2^{1-\gamma})}\left((3-2^{1-\gamma})\sigma-1\right)+\frac{h^2(1-\varepsilon)}{4c_2}\right)\leq
0.
$$
Consequently, from (\ref{ur60}), we have (\ref{ur58.55}) with
$M=M_5(\varepsilon)$.

By multiplying inequality (\ref{ur58.55}) by $\tau$ and by summing
with respect to $s$ from $0$ to $n$, we obtain the a priori estimate
(\ref{ur51}).

{\bf 2)} If $\beta\alpha^{-1}-1<0$ and $\alpha^2-1<0$, then the
second case follows from the first one. Indeed, if we introduce the
notation $y(x,t)=v(1-x,t)$, then the function $v(x,t)$ satisfies
problem ({\ref{ur48}})--(\ref{ur50}).

The solution $v(x,t)$ of problem (\ref{ur48})--(\ref{ur50})
 satisfies the estimate (\ref{ur51}). Consequently, by virtue of
the relations $\|v\|_0^2=\|y\|_0^2$, $\|pv\|_0^2=\|p(1-x)y\|_0^2$
and $v_0^2=\alpha^{-2}y_0^2$, the a priori estimate (\ref{ur33})
holds for the solution of problem (\ref{ur4})--(\ref{ur6}). The
proof of Theorem 5 is complete.

The resulting a priori estimates imply the convergence of the
solution of the difference scheme (\ref{ur4})--(\ref{ur6}) to the
solution of the differential problem (\ref{ur1})--(\ref{ur3}).

The a priori estimates (\ref{ur15}) and (\ref{ur51}) can be obtained
for the case of $\beta=\alpha\neq1$ as well. This follows from the
inequalities
$$
u^2(1,t)=\left(\frac{1}{1-\alpha}\int\limits_{0}^{1}u_x(x,t)dx\right)^2\leq\frac{1}{(1-\alpha)^2}\|u_x\|_0^2.
$$
$$
(y^{(\sigma)}_{N})^2=\left(\frac{1}{1-\alpha}\sum\limits_{i=1}^{N}y^{(\sigma)}_{\bar
x,i}h\right)^2\leq\frac{1}{c_1(1-\alpha)^2}\|\sqrt{a}y^{(\sigma)}_{\bar
x}]|_0^2.
$$

\section{Numerical Results}

Numerical calculations are performed for a test problem when the
function

$$u(x,t)=\left((1-3\alpha)x^3+\alpha x^2+ \alpha x+\alpha\right)\left(t^3-t^2+t+1\right)$$

is the exact solution of the problem (\ref{ur1})--(\ref{ur3}) with
the coefficient $k(x)=e^x$.

The errors ($z=y-u$) and convergence order (CO) in the norms
$|[\cdot]|_0$ and $\|\cdot\|_{C(\bar\omega_{h\tau})}$ at $\sigma=1$
are given in tables 1--7.

Each of tables 1--7 shows that when we take $h^2=\tau^{2-\gamma}$,
as the number of spatial subintervals/time steps is decreased, a
reduction in the maximum error takes place, as expected and the
convergence order of the approximate scheme is $O(h^2)$, where the
convergence order is given by the formula:
CO$=\log_{\frac{h_1}{h_2}}{\frac{\|z_1\|}{\|z_2\|}}$.

Table 8 shows that if $\alpha$ and $\beta$ do not satisfy the
conditions $|\alpha|$, $|\beta| \leq 1$ or $|\alpha|$, $|\beta| \geq
1$, then the difference scheme (\ref{ur4})--(\ref{ur6}) may be
unstable.

\begin{tabular}{lc}
{\bf Table 1}\\
$\gamma=0.5$,  $\alpha=3$, $\beta=2$, $T=1$, $\sigma=1$, $h^2=\tau^{2-\gamma}$\\
 \hline
 $h$ \hspace{12mm}{$\max\limits_{0\leq n\leq N_T}|[z^n]|_0$} \hspace{10mm}{CO in $|[\cdot]|_0$} \hspace{8mm}{$\|z\|_{C(\bar \omega_{h \tau})}$} \hspace{8mm}{CO in $||\cdot||_{C(\bar \omega_{h \tau})}$} \\
\hline
 1/20 \hspace{5mm} $3.03169\cdot10^{-2}$ \hspace{10mm}          \hspace{22mm} $5.50676\cdot10^{-2}$    \hspace{10mm}         \\
 1/40 \hspace{5mm} $7.61510\cdot10^{-3}$ \hspace{10mm}  1.993   \hspace{11mm} $1.38318\cdot10^{-2}$    \hspace{10mm}   1.993 \\
 1/80 \hspace{5mm} $1.90780\cdot10^{-3}$ \hspace{10mm}  1.997   \hspace{11mm} $3.46463\cdot10^{-3}$    \hspace{10mm}   1.997 \\
 \hline
\end{tabular}

\vspace{5mm}

\begin{tabular}{lc}
{\bf Table 2}\\
$\gamma=0.5$,  $\alpha=2$, $\beta=5$, $T=1$, $\sigma=1$, $h^2=\tau^{2-\gamma}$\\
\hline
 $h$ \hspace{12mm}{$\max\limits_{0\leq n\leq N_T}|[z^n]|_0$} \hspace{10mm}{CO in $|[\cdot]|_0$} \hspace{8mm}{$\|z\|_{C(\bar \omega_{h \tau})}$} \hspace{8mm}{CO $||\cdot||_{C(\bar \omega_{h \tau})}$} \\
\hline
 1/20 \hspace{5mm} $6.35368\cdot10^{-3}$ \hspace{10mm}          \hspace{22mm} $7.31523\cdot10^{-3}$    \hspace{10mm}         \\
 1/40 \hspace{5mm} $1.56940\cdot10^{-3}$ \hspace{10mm}  2.017   \hspace{11mm} $1.80908\cdot10^{-3}$    \hspace{10mm}   2.016 \\
 1/80 \hspace{5mm} $3.90276\cdot10^{-4}$ \hspace{10mm}  2.008   \hspace{11mm} $4.49971\cdot10^{-4}$    \hspace{10mm}   2.007 \\
 \hline
\end{tabular}

\vspace{5mm}

\begin{tabular}{lc}
{\bf Tabel 3}\\
$\gamma=0.5$, $\alpha=0.7$, $\beta=0.1$, $T=1$, $\sigma=1$, $h^2=\tau^{2-\gamma}$\\
\hline
 $h$ \hspace{12mm}{$\max\limits_{0\leq n\leq N_T}|[z^n]|_0$} \hspace{10mm}{CO in $|[\cdot]|_0$} \hspace{8mm}{$\|z\|_{C(\bar \omega_{h \tau})}$} \hspace{8mm}{CO in $||\cdot||_{C(\bar \omega_{h \tau})}$} \\
\hline
 1/20 \hspace{5mm} $2.19544\cdot10^{-2}$ \hspace{10mm}          \hspace{22mm} $2.67764\cdot10^{-2}$    \hspace{10mm}         \\
 1/40 \hspace{5mm} $5.50422\cdot10^{-3}$ \hspace{10mm}  1.996   \hspace{11mm} $6.71201\cdot10^{-3}$    \hspace{10mm}   1.996 \\
 1/80 \hspace{5mm} $1.37776\cdot10^{-3}$ \hspace{10mm}  1.998   \hspace{11mm} $1.67992\cdot10^{-3}$    \hspace{10mm}   1.998 \\
 \hline
\end{tabular}

\vspace{5mm}

\begin{tabular}{lc}
{\bf Tabel 4}\\
$\gamma=0.2$, $\alpha=1.1$, $\beta=1.1$, $T=1$, $\sigma=1$, $h^2=\tau^{2-\gamma}$\\
\hline
 $h$ \hspace{12mm}{$\max\limits_{0\leq n\leq N_T}|[z^n]|_0$} \hspace{10mm}{CO in $|[\cdot]|_0$} \hspace{8mm}{$\|z\|_{C(\bar \omega_{h \tau})}$} \hspace{8mm}{CO in $||\cdot||_{C(\bar \omega_{h \tau})}$} \\
\hline
 1/20 \hspace{5mm} $3.85126\cdot10^{-2}$ \hspace{10mm}          \hspace{22mm} $4.38852\cdot10^{-2}$    \hspace{10mm}         \\
 1/40 \hspace{5mm} $9.65615\cdot10^{-3}$ \hspace{10mm}  1.995   \hspace{11mm} $1.10031\cdot10^{-2}$    \hspace{10mm}   1.996 \\
 1/80 \hspace{5mm} $2.42041\cdot10^{-3}$ \hspace{10mm}  1.996   \hspace{11mm} $2.75763\cdot10^{-3}$    \hspace{10mm}   1.996 \\
 \hline
\end{tabular}

\vspace{5mm}

\begin{tabular}{lc}
{\bf Tabel 5}\\
$\gamma=0.2$, $\alpha=0.9$, $\beta=0.9$, $T=1$, $\sigma=1$, $h^2=\tau^{2-\gamma}$\\
\hline
 $h$ \hspace{12mm}{$\max\limits_{0\leq n\leq N_T}|[z^n]|_0$} \hspace{10mm}{CO in $|[\cdot]|_0$} \hspace{8mm}{$\|z\|_{C(\bar \omega_{h \tau})}$} \hspace{8mm}{CO in $||\cdot||_{C(\bar \omega_{h \tau})}$} \\
\hline
 1/20 \hspace{5mm} $3.26779\cdot10^{-2}$ \hspace{10mm}          \hspace{22mm} $3.66507\cdot10^{-2}$    \hspace{10mm}         \\
 1/40 \hspace{5mm} $8.19304\cdot10^{-3}$ \hspace{10mm}  1.996   \hspace{11mm} $9.18862\cdot10^{-2}$    \hspace{10mm}   1.996 \\
 1/80 \hspace{5mm} $2.05366\cdot10^{-3}$ \hspace{10mm}  1.996   \hspace{11mm} $2.30287\cdot10^{-3}$    \hspace{10mm}   1.996 \\
 \hline
\end{tabular}

\vspace{5mm}

\begin{tabular}{lc}
{\bf Tabel 6}\\
$\gamma=0.8$, $\alpha=200$, $\beta=100$, $T=1$, $\sigma=1$, $h^2=\tau^{2-\gamma}$\\
\hline
 $h$ \hspace{12mm}{$\max\limits_{0\leq n\leq N_T}|[z^n]|_0$} \hspace{10mm}{CO in $|[\cdot]|_0$} \hspace{8mm}{$\|z\|_{C(\bar \omega_{h \tau})}$} \hspace{8mm}{CO in $||\cdot||_{C(\bar \omega_{h \tau})}$} \\
\hline
 1/20 \hspace{5mm} $1.27484\cdot10^{0}$ \hspace{10mm}           \hspace{24mm} $2.14188\cdot10^{0}$    \hspace{10mm}         \\
 1/40 \hspace{5mm} $3.18346\cdot10^{-1}$ \hspace{10mm}  2.002   \hspace{11mm} $5.35201\cdot10^{-1}$    \hspace{10mm}   2.001 \\
 1/80 \hspace{5mm} $7.95685\cdot10^{-2}$ \hspace{10mm}  2.000   \hspace{11mm} $1.33790\cdot10^{-1}$    \hspace{10mm}   2.000 \\
 \hline
\end{tabular}

\vspace{5mm}

\begin{tabular}{lc}
{\bf Tabel 7}\\
$\gamma=0.8$, $\alpha=100$, $\beta=200$, $T=1$, $\sigma=1$, $h^2=\tau^{2-\gamma}$\\
\hline
 $h$ \hspace{12mm}{$\max\limits_{0\leq n\leq N_T}|[z^n]|_0$} \hspace{10mm}{CO in $|[\cdot]|_0$} \hspace{8mm}{$\|z\|_{C(\bar \omega_{h \tau})}$} \hspace{8mm}{CO in $||\cdot||_{C(\bar \omega_{h \tau})}$} \\
\hline
 1/20 \hspace{5mm} $6.49129\cdot10^{-1}$ \hspace{10mm}          \hspace{22mm} $1.09160\cdot10^{0}$    \hspace{10mm}         \\
 1/40 \hspace{5mm} $1.62100\cdot10^{-1}$ \hspace{10mm}  2.002   \hspace{11mm} $2.72769\cdot10^{-1}$    \hspace{10mm}   2.001 \\
 1/80 \hspace{5mm} $4.05159\cdot10^{-2}$ \hspace{10mm}  2.000   \hspace{11mm} $6.81875\cdot10^{-2}$    \hspace{10mm}   2.000 \\
 \hline
\end{tabular}

\vspace{5mm}

\begin{tabular}{lc}
{\bf Tabel 8}\\
$\gamma=0.4$, $\alpha=0.1$, $\beta=10$, $T=1$, $\sigma=1$, $h^2=\tau^{2-\gamma}$\\
\hline
 $h$ \hspace{34mm}{$\max\limits_{0\leq n\leq N_T}|[z^n]|_0$} \hspace{44mm} {$\|z\|_{C(\bar \omega_{h \tau})}$} \\
\hline
 1/20 \hspace{27mm} $2.41006\cdot10^{-3}$ \hspace{12mm}          \hspace{22mm} $4.48421\cdot10^{-3}$            \\
 1/40 \hspace{27mm} $5.36386\cdot10^{34}$ \hspace{12mm}          \hspace{22mm} $1.02862\cdot10^{35}$           \\
 1/80 \hspace{27mm} $5.21782\cdot10^{119}$ \hspace{12mm}         \hspace{21mm} $1.0008\cdot10^{120}$           \\
 \hline
\end{tabular}

\section{Conclusion}
The results obtained in the present paper allow to apply the method
of energy inequalities to finding a priory estimate for nonlocal
boundary value problems for the time-fractional diffusion equation
in differential and difference settings. It is interesting to note
that the condition
$$
\sigma\geq\frac{1}{3-2^{1-\gamma}}-\frac{h^2(2-2^{1-\gamma})}{2c_2\tau^\gamma(3-2^{1-\gamma})\Gamma(2-\gamma)}
$$
at $\gamma=1$ turns into the well known condition
$$
\sigma\geq\frac{1}{2}-\frac{h^2}{4c_2\tau}
$$
of the stability of the difference schemes for the classical
diffusion equation.








\end{document}